\documentclass[10pt]{amsart}
\usepackage{amsmath, amsfonts, amsthm, amssymb, verbatim,dsfont}
\usepackage{graphicx}
\usepackage[mathcal]{euscript}
\usepackage[applemac]{inputenc}
\usepackage{dsfont} 
\usepackage{amssymb,mathrsfs}
\usepackage[usenames]{color}
\theoremstyle{plain}
        \newtheorem{theorem}{Theorem}[section]
        \newtheorem{proposition}[theorem]{Proposition}

\theoremstyle{definition}

\theoremstyle{plain}
        
\numberwithin{equation}{section}

\usepackage{amsmath}
\usepackage{verbatim}
\newcommand \be           {\begin{equation}}
\newcommand \ee            {\end{equation}}

\newcommand \RR           {\mathbb{R}}

\newcommand \ZZ           {\mathbb{Z}}

\newcommand \Pbold           {\mathbf{P}} 
\newcommand \PP \Pbold

\newcommand \del           \partial
\newcommand \eps            \epsilon

\newcommand \loc        {{\mathrm{loc}}}

\DeclareMathOperator    \Lip {Lip}
\DeclareMathOperator    \TV  {TV} 
\newcommand{\unj}{u^n_j}

\newcommand{\unnj}{u^{n+1}_j}

\DeclareMathOperator \supp {supp}

\def\XXint#1#2#3{{\setbox0=\hbox{$#1{#2#3}{\int}$}
\vcenter{\hbox{$#2#3$}}\kern-.5\wd0}}

\newcommand{\one}{\mathds{1}}
\let\oldmarginpar\marginpar
\renewcommand\marginpar[1]{\-\oldmarginpar[\raggedleft\footnotesize #1]%
{\raggedright\footnotesize #1}}
\def\build#1_#2^#3{\mathrel{
\mathop{\kern 0pt#1}\limits_{#2}^{#3}}}
\allowdisplaybreaks[1]

\begin{document}

\title[A nonlocal conservation law]{On a nonlocal hyperbolic conservation law arising from a gradient constraint problem}

%

\author{
{ Paulo Amorim}}
\address{
Centro de Matem\'atica e Aplica\c c\~oes
Fundamentais\\
Departamento de Matem\'atica\\
Faculdade de Ci\^encias da Universidade de Lisboa\\
Av. Prof. Gama Pinto 2\\
1649-003 Lisboa, Portugal. 
}
\thanks{Email: pamorim@ptmat.fc.ul.pt. web page: \texttt{http://ptmat.fc.ul.pt/$\sim$pamorim/}
\\Centro de Matem\'atica e Aplica\c c\~oes
Fundamentais,
Departamento de Matem\'atica,
Faculdade de Ci\^encias da Universidade de Lisboa,
Av. Prof. Gama Pinto 2,
1649-003 Lisboa, Portugal. 
}

\begin{abstract} 
In some models involving nonlinear conservation laws, physical mechanisms exist which prevent the 
formation of shocks. This gives rise to conservation laws with a constraint on the gradient of the solution.
We approach this problem by studying a related conservation law with a spatial nonlocal term.
We prove existence, uniqueness and stability of solution of the Cauchy problem for this nonlocal conservation law.
In turn, this allows us to provide a notion of solution to the conservation law with a gradient constraint.
The proof of existence is based on a time-stepping technique, 
and an $L^1$-contraction estimate follows
from stability results of Karlsen and Risebro.
\end{abstract}
 
\maketitle 
 

%
%

\section{Introduction}

We consider the initial value problem for the scalar conservation law with an integral term

\be
\label{10}
\aligned
&\del_t w(t,x) + \del_x \Big( f' \Big(\int_{-\infty}^x w(t, z) dz\Big) g(w(t,x)) \Big) = 0,
\\
& w(0,x) = w_0(x), \qquad x\in \RR, \quad t\ge 0,
\endaligned
\ee
with $f,g : \RR\to \RR$ given functions. The formulation \eqref{10} is motivated by the following problem: to find $u(x,t)$ verifying, in some 
appropriate sense, the following \emph{conservation law
with a gradient constraint,}
\be
\label{20}
\aligned
&\del_t u  + \del_x (f(u)) =0, \qquad u(0,x) = u_0(x),
\\
& | \del_x u(t,x)| \le M, \qquad x\in \RR, \quad t\ge 0,
\endaligned
\ee
for some $M>0$.

Problem \eqref{20} has no ready interpretation in the scope of conservation laws. Indeed, as is
well known, a nonlinear conservation law will, even for smooth initial data, develop discontinuities in
finite time, whose onset is preceded by a blowup of the spatial derivative. Therefore, it is hopeless to
seek solutions verifying problem \eqref{20}, without some additional information. 

We now describe our interpretation of problem \eqref{20}, which leads to the formulation \eqref{10}. 
Setting $w := \del_x u$, then spatial differentiation of \eqref{20} gives formally
\be
\label{15}
\aligned
&\del_t w + \del_x \Big( f' \Big(\int_{-\infty}^x w(t, z) dz\Big) w \Big) = 0,
\\
& |w(t,x)| \le M,
\endaligned
\ee
so that the constraint on the gradient of $u$ now acts on the function $w$ itself. We provide a formulation for this problem
by replacing the term $w$ appearing inside the spatial derivative in \eqref{15} by a function \emph{with compact support} $g(w)$.
This, as we shall see, limits the growth of $|w|$ (and, formally, of the gradient $|\del _x u|$). Thus, we arrive at problem \eqref{10}.
Note that the desired constraint which appears in \eqref{15} is now ensured by the compactness of the support of the function $g(w)$ in \eqref{10}.

The motivation for considering problems \eqref{10},\eqref{20} comes from the fact that in some situations involving
conservation laws, there may exist some physical mechanism preventing the formation of shocks. For instance, models of 
granular motion \cite{AS} and superconductors \cite{Prigozhin} have such characteristics.

Other approaches to this type of problem include interpreting \eqref{20} as a quasi-variational problem, which can be shown to be
well approximated by a viscous conservation law with a highly singular viscosity function. This approach was 
carried out successfully by Rodrigues and Santos \cite{RS} in the context of free boundary
problems. It would be interesting to compare the solutions obtained in the present work
to the ones in that paper. Let us also mention the work of L\'evi \cite{Levi1}, in which a conservation law
with a positivity constraint arising in oil reservoir dynamics is solved by a penalization method.

Even though these approaches may be fruitful, our standpoint is different and simpler. 
We find it of interest to consider the hyperbolic formulation
\eqref{10}, so that the problem may be treated using the techniques of hyperbolic conservation laws.
Moreover, equations of the form \eqref{10} are of interest in themselves, regardless of the motivation
given here. Indeed, it is now well-known that a variety of physical phenomena are better described
through the introduction of integral terms. For example, we mention the recent work of Amadori and Shen \cite{AS}, 
where the authors consider a problem arising in the study of granular motion.
The present work shares some of its
techniques with \cite{AS}, namely the time discretization used to deal with the integral term.

Another domain of application of conservation laws with nonlocal term is given by the modeling of pedestrian flows (see \cite{CHM} and the references 
therein). In that context, the equation \eqref{15} with a nonlocal term models the density of pedestrians evolving in time, where the nonlocal 
velocity function $f'(\int_{-\infty}^x w)$ translates the fact that pedestrians adjust their speed according to a perceived average of the density on some domain,
in this case $(-\infty,x)$.

Let us also briefly mention that nonlocal in time terms are also relevant in hyperbolic conservation laws; in this respect, 
we only mention the pioneering work of Dafermos \cite{Dafermos}.

In this work, we provide a well-posedness theory for problem \eqref{10}. Existence follows from a time-stepping technique also employed 
by Amadori and Shen in 
\cite{AS}. We point out that in that paper, even though the authors consider a related problem, the assumptions therein are very different from ours.
The practical consequence is that in \cite{AS} the authors may rely on previous works in which the important auxiliary problem \eqref{200} below 
can be treated using the decomposition into a $2\times2$ system of conservation laws (see, for instance, Klausen and Risebro \cite{KlausenRisebro}).
This framework is by now well-established but contains assumptions which, although convenient in \cite{AS}, do not suit our framework.
Thus, we need to rely instead on more recent results of Karlsen and Risebro \cite{KR1,KR2} and Chen and Karlsen \cite{CK}. Based on these results,
we are able to prove an $L^1$-stability property for the solutions.

Relying on the well-posedness framework developed for equation \eqref{10}, we are able, in Theorem~\ref{thm17} below, to give a precise meaning
to the solution of problem \eqref{20} with a gradient constraint.

\subsection{Main results and assumptions}
We now present our main results. First, we state precisely the assumptions on the functions appearing in \eqref{10}, 
\eqref{20}. Henceforth, we assume
\be
\label{100}
f : \RR \to \RR \quad \text{is a $C^3$ function with }  \| f'' \|_{L^\infty(\RR)} \le C,
\ee
\be
\label{110}
g : \RR \to \RR \in \Lip(\RR), \quad \supp g \subset [-M, M], \quad M>0,
\ee
\be
\label{120}
w_0  : \RR \to \RR \in  L^1(\RR) \cap BV(\RR), \quad |w_0| \le M.
\ee

It is worth noting that the property \eqref{110} of boundedness of the support of $g$ is not a drawback, but rather a
feature, of this work. Indeed, without such an assumption --- if we had, for instance, $g(w) =w$ --- then 
\eqref{10} would be the equation obtained by spatial differentiation of the conservation law \eqref{20} (without any constraints).
Thus, $w$ would naturally blow up in the $L^\infty$ norm in finite time, being the derivative of a function developing discontinuities.

Therefore, as we shall see, condition~\eqref{110} is the cornerstone of our attempt at interpreting problem \eqref{20}. 
One natural choice of function $g$ is thus a smoothed version of the truncated identity function; in that case,
whenever $w$ remains inside the support of $g$, then by integration, we may recover $u$ from $w$, and $u$ solves
\eqref{20}. Since we will prove below that $|w| \le M$ for all $(t,x)$, it is reasonable to give \eqref{10} as
an interpretation of problem~\eqref{20}. (Of course here, the function $g$ may have some more general expression, as long 
as it satisfies condition~\eqref{110}). This outline is made more precise in Theorem~\ref{thm17} below.

Let $f,g$ and $w_0$ verify the assumptions \eqref{100}--\eqref{120}. We say that a function $w \in L^\infty([0,T]; L^1\cap BV(\RR))$ is an entropy solution to the problem \eqref{10} if
it is a weak entropy solution, in the standard sense, of the Cauchy problem
\be
\label{}
\aligned
&\del_t w + \del_x \big( k(t,x)  g(w)\big) = 0,
\\
&w(0,x) = w_0(x),  \qquad x\in \RR, \quad t\in [0,T],
\endaligned
\ee
with $k(t,x)$ given by 
\[
\aligned
k(t,x) = f'\Big( \int_{-\infty}^x w(t,z) \,dz \Big).
\endaligned
\]

Our first result establishes existence of a solution to problem~\eqref{10}.
\begin{theorem}
\label{thm10}
Let the functions $f$, $g$ and $w_0$ satisfy the assumptions \eqref{100}--\eqref{120}, and let $T>0$. Then, there exists
an entropy solution $w(t,x)$ to the nonlocal conservation law \eqref{10}. This solution 
satisfies, for all $0\le t\le T$, the estimates
\be
\label{40}
\aligned
\|w\|_{L^\infty((0,T)\times \RR)} \le M,
\endaligned
\ee
\be
\label{60}
\aligned
\|w(t)\|_{L^1(\RR)} \le \|w_0\|_{L^1(\RR)},
\endaligned
\ee
and
\be
\label{65}
\aligned
\TV(w(t)) \le C(t),
\endaligned
\ee
for some continuous function $C(t)$.

\end{theorem}

We also prove the following $ L^1$-stability result:
\begin{theorem}
\label{thm15}
Let the functions $f$, $g$ and $w_0$ satisfy the assumptions \eqref{100}--\eqref{120}. Then,
the following $L^1$-stability property is valid: if $w,v$ are entropy solutions which verify the uniform estimates \eqref{40},\eqref{60}, with
initial data $w_0,v_0$ satisfying \eqref{120}, then for all $t>0$ we have
\be
\label{50}
\aligned
\| w(t) - v(t) \|_{L^1(\RR)} \le e^{Ct}\| w_0 - v_0 \|_{L^1(\RR)},
\endaligned
\ee
where $C$ is a constant depending only on $g$, $f$ and the initial data.
\end{theorem}
In particular, Theorem~\ref{thm15} gives uniqueness of entropy solution to \eqref{10} within the class of functions satisfying \eqref{40},\eqref{60}.

We postpone the proofs of Theorems \ref{thm10} and \ref{thm15} to later sections.

Let us now translate the results of Theorem \ref{thm10} into a result for the conservation law with gradient constraint \eqref{20}.

\begin{theorem}
\label{thm17}
Suppose that the function $g$ in \eqref{110} has the form $g(w) = w h_\eps(w)$, where $h = h_\eps$ is some regularization of the characteristic function
$\one_{[-M,M]}$, with $h=1$ on $[-M+\eps,M-\eps]$, and let $w$ be the solution of \eqref{10} given by Theorems \ref{thm10} and \ref{thm15}. 
Define the following sets contained in $(0,T) \times \RR$,
\[
\aligned
&I_\eps = \{ (t,x) : |w| \le M-\eps \},
\\
&J_\eps = \{ (t,x) : |w| = M \},
\\
&K_\eps = ((0,T) \times \RR) \setminus (I_\eps \cup J_\eps).
\endaligned
\]

Then, the function 
$u := \int_{-\infty}^x w(t,y) \,dy$ solves the conservation law with gradient constraint \eqref{20} in the sense that
$u$ is an entropy solution of
\[
\aligned
&\del_t u  + \del_x (f(u)) =0,
\\
& | \del_x u | < M, 
\endaligned
\qquad \text{ on}\quad I_\eps,
\]
$u$ verifies
\[
\aligned
| \del_x u| = M \qquad \text{ on}\quad J_\eps,
\endaligned
\]
and $u$ solves
\[
\aligned
&\del_t u  + \del_x (f(u))h(w) =0 
\endaligned
\]
on the transition layer $K_\eps$.
Furthermore, $u$ verifies the estimate
\be
\label{52}
\aligned
u \in L^\infty(0,T; W^{1,\infty}(\RR)).
\endaligned
\ee
In particular, $u$ is continuous on $\RR$ for each $t$.
\end{theorem}
Theorem \ref{thm17} is an easy consequence of Theorem \ref{thm10} and of arguments similar to \cite[Lemma 2.2.1]{BouchutJames}, 
and so we omit a detailed proof. In particular, the estimate
\eqref{52} follows from \eqref{40} and \eqref{60}. Note also that the equation $\del_t u  + \del_x (f(u))h(w) =0$ verified by $u$ on the
transition layer $K_\eps$ is obtained from \eqref{10} when $g = w h$.

\subsection{Comments and remarks}
\begin{enumerate}
\item The assumption in Theorem~\ref{thm15} that the solutions verify condition \eqref{60} can easily be relaxed 
to $\|w,v\|_{ L^1(\RR)} \le C(t),$ and the result of Theorem~\ref{thm15} (and so, in particular, the estimate \eqref{60} itself) remain valid. 

\item The problem \eqref{20} with a gradient constraint and its solution $u$ given by Theorem~\ref{thm17} can be seen from 
the viewpoint of free boundary problems, as in \cite{RS}. Indeed, one can view the sets $I_\eps$ and $J_\eps$ in the statement of
Theorem~\ref{thm17} as two domains separated by a thin \emph{transition layer} $K_\eps$, which should become a free boundary 
as $\eps$ tends to zero. In $I_\eps$, $u$ solves the conservation law; in $J_\eps$, $u$ solves the Hamilton--Jacobi equation
$| \del_x u| =M$.

\item Since the function $g$ in \eqref{10} is required to be smooth, we cannot simply take
a (discontinuous) truncation of the identity function on $[-M,M]$. Therefore, our result depends on some small
smoothing parameter used to regularize the identity function on $[-M,M]$. On the other hand, we allow for more general smooth 
functions $g$ supported in $[-M,M]$, if we do not wish to see \eqref{10} as an approximation of the constrained problem \eqref{20}. 
In this direction, it would be interesting to study the applicability to this problem of the results in
\cite{Checos} or \cite{DF} dealing with conservation laws with a discontinuous flux function.

\item It can be seen from a careful analysis of the proof of Theorem~\ref{thm10} that, in the case where $f$ is convex, it is only necessary 
to suppose that the support of $g$ is bounded from below in order to obtain the $ L^\infty$ bound \eqref{40}.

\item Several extensions of the results in this paper are possible, namely,
the extension to the (more realistic) situation where the gradient constraint depends on $t$ and $x$. Also,
the numerical treatment of \eqref{10} would be interesting to study. This would shed light on the
relation between the solutions obtained in the present work and the ones obtained in \cite{RS}. We plan to address these 
questions in further work.
\end{enumerate}

\section{Preliminary results}
In this section, we establish some auxiliary results.
We now consider problems of the type
\be
\label{200}
\aligned
\del_t w + \del_x (k(t,x) g(w)) = 0,
\endaligned
\ee
where now $k$ is a \emph{fixed} function which we assume has the following regularity:
\be
\label{230}
\aligned
k(t, \cdot) \in W^{1,1}\cap W^{1,\infty} \cap L^\infty (\RR), \qquad \del_x k(t, \cdot) \in BV(\RR), \quad 
\text{uniformly in }t.
\endaligned 
\ee
In view of establishing a stability result with respect to the function $k$, let us also introduce
the similar problem
\be
\label{210}
\aligned
\del_t v + \del_x (l(t,x) g(v)) = 0,
\endaligned
\ee
with the function $l(t,x)$ verifying the same assumptions as $k$.

We recall that the problem \eqref{200} with the assumption \eqref{230} is not contained in the classical 
work of Kruzhkov \cite{K}. To the author's knowledge, the first well-posedness results with this kind of rough coefficients are to be found in the works of Karlsen and Risebro \cite{KR1,KR2}.

The following result is a slight extension of \cite[Theorem 1.1]{KR1}, \cite[Theorem~1.3]{KR2} and \cite[Theorem~6.1]{CK}.
Thus, we shall only point out
where the proof differs from the ones in those papers, providing the necessary arguments as needed. 

\begin{theorem}
\label{thm20}
Let $k, g, w_0, v_0$ be functions satisfying assumptions \eqref{110},\eqref{120} and \eqref{230}.
Then, there exists a unique entropy solution to the problem \eqref{200} with initial data $w_0$. 

Suppose now that the functions $k,l$ depend only on the space variable $x$. 
Then, the unique solution $w$ of \eqref{200} 
belongs to $ L^\infty([0,T] ; BV(\RR))$ and
verifies the uniform estimate
\be
\label{250}
\aligned
\|w\|_{L^\infty((0,T)\times \RR)} \le M
\endaligned
\ee
where $M$ is such that $\supp g \subset [-M,M]$.
Moreover,
the following continuous dependence estimate is valid: If $w$ is a solution of problem \eqref{200} with 
initial data $w_0$, and $v$ is a solution of problem \eqref{210} with initial data $v_0$, with $w,v \in  L^\infty([0,T] ; BV(\RR)),$ then for 
$t_1 \le t_2 \le T$, we have
\be
\label{260}
\aligned
&\| w(t_2) - v(t_2) \|_{L^1(\RR)} \le \| w(t_1) - v(t_1) \|_{L^1(\RR)} 
\\
&\qquad +\int_{t_1}^{t_2}  \Lip_g \| k -l \|_{L^\infty(\RR)} \TV(w(\tau))\wedge \TV(v(\tau))  + 
M \TV( k - l)    \, d\tau,
\endaligned
\ee
where $a\wedge b := \min\{a,b\}$.
In particular, any solution $w$
verifies the $L^1$-contraction property for every $t>0$,
\be
\label{255}
\aligned
\|w(t)\|_{L^1(\RR)} \le \|w_0\|_{L^1(\RR)}.
\endaligned
\ee

\end{theorem}
\proof
The existence property in Theorem \ref{thm20} is just Theorem~6.1 of \cite{CK}.
Regarding the remaining assertions of Theorem \ref{thm20},
the only points of difference with the result \cite[Theorem 1.1]{KR1} are the $ L^\infty$ bound \eqref{250}, and the property $w \in L^\infty([0,T]; BV(\RR)),$
which is absent from the statement of that theorem. 

In what follows, we refer the reader to \cite{KR1}, especially Section 3 in that paper for more details. There, the authors consider
a numerical approximation of the equation \eqref{200} with $k = k(x)$ (actually, a more general version of \eqref{200}) which
they use to prove existence of solution. 

First, observe that the discrete version of the property $w \in L^\infty([0,T]; BV(\RR))$ can be found in \cite[p.253]{KR1}. This gives a similar bound
for the exact solution after passing to the limit on the discretization parameter used in that proof. We don't bother to write the
exact estimate since it is not precise enough for our purposes and will be refined later on.

Let us now turn to the estimate \eqref{250}.
It is enough to establish an estimate of the form \eqref{250} for the
approximate solutions employed in \cite{KR1}, obtaining in the limit the corresponding estimate for the exact solution.

In order to use the same notations as \cite{KR1}, let $u^n_j$ denote the solution of a numerical approximation of equation \eqref{200}
associated with the discretization parameters $\Delta t$ and $\Delta x$. Here,
$n$ represents the time level and $j$ the spatial point in some mesh with nodes $\{ x_{j+1/2}\}_{j\in\ZZ}$.
Also, let $U^n = \sup_{j\in\ZZ} |\unj|$. 

From \cite[p.~250]{KR1}, and as a consequence of the monotonicity of the numerical scheme employed, we have
\[
\aligned
|\unnj |  &\le U^n + \lambda| k^n_{j+1/2} - k^n_{j-1/2} | |g(U^n)|,
\endaligned
\]
where $\lambda = \frac{\Delta t }{\Delta x}$ and $k^n_{j+1/2}$ denotes some discretization of the function $k(t,x)$ verifying 
$| k^n_{j+1/2} - k^n_{j-1/2} | \le \Delta x\Lip_x k.$ 

Let $n$ be the first time level for which $U^{n+1} > M$, where $\supp g \subset [-M,M]$. Then, $g(U^{n+1}) = 0$, so we can write
\[
\aligned
U^{n+1} &\le U^n + \lambda \sup_{j\in\ZZ}| k^n_{j+1/2} - k^n_{j-1/2} | |g(U^n) - g(U^{n+1})|,
\endaligned
\]
which gives
\[
\aligned
U^{n+1} - U^n \le  \lambda \sup_{j\in\ZZ} | k^n_{j+1/2} - k^n_{j-1/2} | (U^{n+1} - U^n ) \Lip g .
\endaligned
\]
Using the CFL condition $\lambda \Lip g \le 1$ (which we can use instead of \cite[(3.8)]{KR1} due to the particular form of the flux term,
$f(k,u) = k(x)u$), we obtain
\[
\aligned
1 \le \sup_{j\in\ZZ} | k^n_{j+1/2} - k^n_{j-1/2} | \le \Delta x \Lip_x k,
\endaligned
\]
which is seen to be a contradiction by, say, taking $\Delta x$ sufficiently small so that $\Delta x( \Lip_x k +1) \le 1$. Thus, $U^n \le M$ for all $n$, which is the
discrete version of \eqref{250}.

We now turn to the estimate \eqref{260}. This is simply a rewriting of \cite[Theorem~1.3]{KR2}. Remark that by carefully analyzing the proof, in particular \cite[p.1011]{KR2},
we can see that this result may indeed be formulated with integration in time instead of multiplication by $t$.

\endproof

\section{Existence of solution to the nonlocal conservation law}
\label{sec30}

\subsection{A time-stepping technique}
We will now consider the problem \eqref{10}. Our strategy, inspired by Amadori and Shen 
\cite{AS}, is to consider a time-stepping technique to obtain a solution $w$ of \eqref{10} by
compactness of a family of approximate solutions $w^\delta(t,x)$. The idea is that for each fixed $T,\delta>0$,
we set $t_n:= n\delta$, $n=1,2,\dots$, and define a function $w^{\delta,1}(t,x)$, in the interval $t \in [0, t_1)$ 
as the unique entropy solution of the problem
\[
\aligned
&\del_t w^{\delta,1} + \del_x \Big( f'\Big(\int_{-\infty}^x w_0(y) dy \Big) g(w^{\delta,1}) \Big) = 0,
\\
& w^{\delta,1}(0,x) = w_0(x)
\endaligned
\]
given by Theorem \ref{thm20}. Next, we set $w^{\delta,2}(t,x)$ in the interval $t \in [t_1, t_2)$ as the 
unique entropy solution of the problem
\[
\aligned
&\del_t w^{\delta,2} + \del_x \Big( f'\Big(\int_{-\infty}^x w^{\delta,1}(t_1-,y) dy \Big) g(w^{\delta,2}) \Big) = 0,
\\
& w^{\delta,2}(t_1,x) = w^{\delta,1}(t_1-,x),
\endaligned
\]
and so on, so that $w^{\delta,n+1}(t,x)$ is defined in the interval $t\in [t_{n}, t_{n+1})$ and is a solution of
\[
\aligned
&\del_t w^{\delta,n+1} + \del_x \Big( f'\Big(\int_{-\infty}^x w^{\delta,n}(t_n-,y) dy \Big) g(w^{\delta,n+1}) \Big) = 0,
\\
& w^{\delta,n+1}(t_n,x) = w^{\delta,n}(t_n-,x).
\endaligned
\]
Finally, we set
\be
\label{300}
\aligned
w^\delta(t,x) = \sum_{i=1}^n \chi_{[t_n, t_{n+1}]} w^{\delta,n}(t,x).
\endaligned
\ee

\subsection{Estimates for the time-stepping approximate solution}
Next, we shall prove the crucial estimates on $w^\delta$.
\begin{proposition}
\label{prop30}
Let $ w^\delta (t,x)$ be defined by \eqref{300}. Then, for each $t\in[0,T]$, we have the estimates
\be
\label{360}
\aligned
\| w^\delta \|_{L^\infty((0,T)\times\RR)} \le M,
\endaligned
\ee
\be
\label{370}
\aligned
\| w^\delta(t) \|_{L^1(\RR)} \le \| w_0 \|_{L^1(\RR)}.
\endaligned
\ee
\be
\label{350}
\aligned
\TV (w^\delta(t)) \le (1 + \TV(w_0)) e^{C t},
\endaligned
\ee
with $ C = C \big(M, \Lip g, \Lip f', \Lip f'', \| w_0\|_{L^1(\RR)} \big).$

\end{proposition}

\proof
First of all, note that the estimates \eqref{360},\eqref{370} are immediate from the corresponding ones
in Theorem~\ref{thm20}. We now prove the total variation estimate \eqref{350}. 
Observe that from Theorem~\ref{thm20}, we know that $ w^\delta$ has bounded variation on each interval $[t_n,t_{n+1})$.
However, we need a precise estimate on each such interval in order to pass to the limit in $\delta$.

For this, let $h>0$
be a small parameter, and set
$\overline w^\delta(x,t) :=  w^\delta (t,x+h)$. Then, it is easy to verify that on each interval
$[t_n,t_{n+1}]$, $\overline w^\delta$ is a solution of
\[
\aligned
&\del_t \overline w^{\delta,n+1} + \del_x \big( \overline {k^n}(t,x)
g(\overline w^{\delta,n+1}) \Big) = 0,
\\
& \overline w^{\delta,n+1}(t_n,x) = \overline w^{\delta,n}(t_n-,x),
\endaligned
\]
where we set 
$$
k^n(t,x) := f'\Big( \int_{-\infty}^{x} w^{\delta,n} (t_n-,y) \, dy \Big),
$$ 
and 
$$
\overline {k^n}(t,x) := k^n(t,x+h) = f' \Big( \int_{-\infty}^{x+h} w^{\delta,n} (t_n-,y) \, dy  \Big).
$$
Note that $k^n$ is constant in time in each interval $[t_n, t_{n+1})$.
For clarity's sake, we omit the superscript $\delta$ from the remainder of this proof. Now, since on each interval $[t_n, t_{n+1})$ the 
functions $k^n$ and $\overline{k^n}$ are constant in time, we may apply the continuous dependence estimate \eqref{260}
to compare $w$ with $\overline w$ on each interval $ [t_n, t_{n+1})$, with $k = k^n$ and $l = \overline{k^n}$.
Note that by \eqref{120}, the expression of $k$, and the estimates in Theorem~\ref{thm15} (applied in previous time steps), the function
$k^n$ indeed verifies the conditions \eqref{230}, so we may apply \eqref{260}.

Thus, it is convenient to first estimate the terms appearing in \eqref{260}. We easily find
\[
\aligned
\|  k^n - \overline{k^n} \|_{L^\infty(\RR)} \le \Lip f' \sup_{x \in\RR} \Big| \int_{x}^{x+h} w( t_n, y)\, dy \Big| \le h \Lip f' \| w (t_n) \|_{L^\infty(\RR)}
\endaligned
\]
and
\[
\aligned
&\hspace{-10pt} \del_x( k^n - \overline{k^n}) 
\\
& \le \Big| f'' \Big( \int_{-\infty}^x w(t_n,y) \, dy \Big) \, w(t_n, x) - f'' \Big( \int_{-\infty}^{x+h} w(t_n,y) \, dy \Big) \, w(t_n, x+h) \Big|
\\
&\le \Lip f'' h \|w(t_n) \|_{L^\infty(\RR)} | w(t_n, x)| + \| f'' \|_{L^\infty(\RR)} | w(t_n, x) -  \overline w (t_n, x) |,
\endaligned
\]
and so, using \eqref{360},\eqref{370},
\[
\aligned
\TV ( k^n - \overline{k^n}) & \equiv \| \del_x( k^n - \overline{k^n}) \|_{L^1(\RR)} 
\\
&\le  h \Lip f''  M \| w_0 \|_{L^1(\RR)} + \| f'' \|_{L^\infty(\RR)} \| w(t_n) - \overline w (t_n) \|_{L^1(\RR)} .
\\
\endaligned
\]
Thus, the continuous dependence estimate \eqref{260} gives
\be
\label{}
\aligned
&\frac1h \| w(t_{n+1}) - \overline w (t_{n+1}) \|_{L^1(\RR)} \le  \frac1h \| w(t_n) - \overline w (t_n) \|_{L^1(\RR)}  
\\
&\qquad+ \Lip g \Lip f' M \int_{t_n}^{t_{n+1}} \TV (w (\tau)) \, d\tau 
\\
&\qquad+ M \int_{t_n}^{t_{n+1}}\Lip f'' M \|w_0 \|_{L^1(\RR)} + \| f''\|_{L^\infty(\RR)} \TV(w(t_n)) \, d\tau.
\endaligned
\ee
Summing in $n = 1, \dots, N$ such that $t = N\delta$, and recalling the definition of the total variation,
\[
\aligned
TV(w) = \lim_{h\to0} \frac1h \int_{\RR} |w(x) - w(x+h)| \,dx,
\endaligned
\]
we find
\[
\aligned
\TV(w(t)) \le \TV(w_0) + C \int_0^t 1 + \sup_{\tau \in (0,t)} \TV(w(\tau)) \, d\tau
\endaligned
\]
for some appropriate constant $C = C\big(M, \Lip g, \Lip f', \Lip f'', \| w_0\|_{L^1(\RR)} \big)$. By Gronwall's lemma
applied to $1 + \sup_{t \in (0,t)} \TV(w(t)) $, this gives
\[
\aligned
\TV(w(t)) \le \big( 1+ \TV(w_0) \big) e^{Ct},
\endaligned
\]
which is \eqref{350}. This completes the proof of Proposition~\ref{prop30}.
\endproof

\subsection{Proof of Theorem \ref{thm10}: Existence of solution to the full problem}
We are now in a position to obtain a solution to the problem \eqref{10}, as a consequence of the
estimates in Proposition \ref{prop30}. We sketch the arguments, since they are quite standard. Indeed, the approximate solutions $w^\delta$ are 
in the space $ L^\infty( [0,T] ; BV(\RR))$, and so by a well-known compactness result, there is a
subsequence (still labelled $w^\delta$) converging in $L^1_\loc([0,T] \times \RR)$ and a.~e.~to a function
$w(t,x)$. 

Now, from the bounds \eqref{360},\eqref{370}, we find 
$$|k^n(t,x)| \le \sup_{[-\| w_0\|_{ L^1(\RR)}, \| w_0\|_{ L^1(\RR)}]} |f'|,
\qquad |g( w^\delta)| \le \sup_{[-M,M]} |g|, $$
and so we may use Lebesgue's dominated convergence theorem
to pass to the limit $\delta\to 0$ on the integral formulation of the conservation law. Thus, $w$ is an entropy solution of
problem \eqref{10}. This completes the proof of Theorem~\ref{thm10}. \qed

\section{Uniqueness of solution to the nonlocal conservation law}

\subsection{Proof of Theorem \ref{thm15}}
The main step in the proof of Theorem \ref{thm15} is to generalize the continuous dependency estimate \eqref{260}
to deal with velocity functions $k,l$ depending on the time $t$ as well as $x$. Thus, let us consider
a solution $w$ to problem \eqref{200} and $v$ a solution to problem \eqref{210} with initial data $w_0$ and $v_0$, respectively. We have
\[
\aligned
\del_t w + \del_x (k(t,x) g(w)) = 0,\qquad \del_t v + \del_x (l(t,x) g(v)) = 0.
\endaligned
\]
Now, as in Section \ref{sec30}, we consider for a fixed $t>0$, and for each $\delta>0$ the approximations $k^\delta(t,x), l^\delta(t,x)$,
constant on each time interval $[t_n, t_{n+1})$ (where $t_n = n\delta$), given by $k^\delta(t,x) = k(t_n,x),$ $l^\delta(t,x) = l(t_n,x).$
According to Theorem~\ref{thm20}, there exists on each $[t_n, t_{n+1})$ a solution $w^n$ of the problem
\[
\aligned
\del_t w^n + \del_x (k^\delta(t,x) g(w^n)) = 0,\qquad w^n(t_n,x) = w^{n-1}(t_n-, x),
\endaligned
\]
and we define the functions $w^\delta, v^\delta$ on $[0,t]$ in the same way as in \eqref{300}. Set
\[
\aligned
\Psi^\delta(t) := \| w^\delta(t) - v^\delta(t) \|_{L^1(\RR)}.
\endaligned
\]
Then, the continuous dependence estimate \eqref{260} gives
\[
\aligned
\Psi^\delta (t_{n+1}) &\le \Psi^\delta (t_n) + \Lip g \| k(t_n) - l(t_n) \|_{L^\infty(\RR)} \int_{t_n}^{t_{n+1}} \TV (w^\delta(\tau)) \wedge \TV (v^\delta(\tau)) \,d\tau
\\
& + \int_0^t M \TV(k(t_n) - l(t_n)) \, d\tau.
\endaligned
\]
summing in $n = 0,\dots ,N$ so that $N\delta = t$ (which, for simplicity, we can assume to be the case), and taking the supremum in time of the quantities
involving $k$ and $l$, we find
\be
\label{450}
\aligned
\Psi^\delta (t) &\le \Psi^\delta(0) + \Lip g  \int_{0}^{t}\sup_{[0,\tau)}\| k - l \|_{L^\infty(\RR)} \TV (w^\delta(\tau)) \wedge \TV (v^\delta(\tau)) \, d\tau 
\\
&+ \int_0^t M \sup_{[0,\tau)} \TV(k(\tau) - l(\tau)) \, d\tau.
\endaligned
\ee
We must now estimate $\TV(w^\delta)$ and $\TV(v^\delta)$. As in the proof of Proposition~\ref{prop30}, we use the continuous dependence estimate
\eqref{260} with $ w^\delta $ and $ \overline w^\delta(t,x) := w^\delta(t,x+h)$. Note that $\TV ( k - \overline k) \le h \TV (\del_x k)$ and that
$\| k - \overline k\|_{ L^\infty(\RR)} \le h \Lip k$,
so we obtain
\[
\aligned
&\frac1h \| w^\delta(t) - \overline w^\delta(t) \|_{ L^1(\RR)} \le \frac1h \| w^\delta(0)- \overline w^\delta(0) \|_{ L^1(\RR)} 
\\
&\quad + \int_0^t \Lip g \TV ( w^\delta(\tau)) \sup_{[0,\tau)} \Lip k + M \sup_{[0,\tau)} \TV(\del_x k) \, d\tau.
\endaligned
\]
Using the definition of the total variation and Gronwall's lemma, we find
\be
\label{460}
\aligned
\TV (w^\delta(t)) \le \big( \TV(w^\delta(0)) + t \sup_{[0,t)}\TV (\del_x k) \big) e^{c_1t}, 
\endaligned
\ee
with $c_1 = \Lip g \sup_{[0,t)} \Lip k. $ A similar estimate is valid for $ \TV(v^\delta(t))$, with $l$ instead of $k$. Plugging these estimates in \eqref{450}
we find
\be
\label{465}
\aligned
\Psi^\delta (t) &\le \Psi^\delta(0) + \Lip g \int_{0}^{t}  \sup_{[0,\tau)}\| k - l \|_{L^\infty(\RR)} \Theta (\tau, w_0, k) \wedge \Theta (\tau, v_0, l)  \, d\tau
\\
&+ \int_0^t M \sup_{[0,\tau)} \TV(k(\tau) - l(\tau)) \, d\tau,
\endaligned
\ee
with $\Theta(\tau, w_0, k)$ given by the right-hand side of \eqref{460}. Recall that $\Psi^\delta = \| w^\delta - v^\delta\|_{ L^1(\RR)}$.

Now, in view of the regularity of $k,l$ (see \eqref{230}), and the estimate \eqref{460}, it is a simple matter to prove that (for a fixed $t>0$) $ w^\delta$ and $ v^\delta$ converge (as $\delta \to 0$) in $ L^1(\RR)$ to $w,v$ which are
the entropy solutions of the problems \eqref{200} and \eqref{210}. Thus, passing to the limit $\delta \to 0$ on the estimate \eqref{465}, we obtain
with $\Psi = \| w - v\|_{ L^1(\RR)}$,
\be
\label{470}
\aligned
\Psi (t) &\le \Psi (0) + \Lip g \int_{0}^{t} \sup_{[0,\tau)}\| k - l \|_{L^\infty(\RR)} \Theta (\tau, w_0, k) \wedge \Theta (\tau, v_0, l)  \, d\tau
\\
&+ \int_0^t M \sup_{[0,\tau)} \TV(k(\tau) - l(\tau)) \, d\tau,
\endaligned
\ee
with $\Theta(\tau, w_0, k)$ given by the right-hand side of \eqref{460}.

With the estimate \eqref{470} in hand, we are now in a position to conclude the proof of Theorem~\ref{thm15}. 
In the remainder of the proof, let us denote by $w$ and $v$ two solutions of the conservation law with integral term \eqref{10}, without causing
confusion with what precedes. Thus, $w,v$ verify
\[
\aligned
&\del_t w + \del_x \Big( f' \Big(\int_{-\infty}^x w(t, z) dz\Big) g(w) \Big) = 0,
\\
& \del_t v + \del_x \Big( f' \Big(\int_{-\infty}^x v(t, z) dz\Big) g(v) \Big) = 0,
\endaligned
\]
with initial data $ w(0,x) = w_0(x),$ $v(0,x) = v_0(x)$. $w$ and $v$ may therefore be seen as solutions of the problems \eqref{200} and \eqref{210}, respectively,
where $k(x,t) = f' \Big(\int_{-\infty}^x w(t, z) dz\Big)$ and $l(x,t) = f' \Big(\int_{-\infty}^x v(t, z) dz\Big)$.

Our goal is to apply the estimate \eqref{470}. For this, it is convenient to estimate the terms appearing in \eqref{470} which involve $k$ and $l$.
Thus,
\[
\aligned
\|k(t,x) - l(t,x) \|_{ L^\infty(\RR)} &\le \Lip f' \int_{-\infty}^x | w(t, y) - v(t,y) | \, dy 
\\
&\le \Lip f' \| w(t) - v(t) \|_{ L^1(\RR)},
\endaligned
\]
\[
\aligned
\TV(\del_x k) &\le \int_\RR \big| f''\Big( \int_{-\infty}^x w(t,y) \, dy \Big) w(t,x) \big| \, dx \le \|f''\|_{ L^\infty(\RR)} \| w_0\|_{ L^1(\RR)},
\endaligned
\]
and
\[
\aligned
\TV( k(t) - l(t)) &\le \int_{\RR}  \big| f'' \Big( \int_{-\infty}^x w(t,y) \, dy \Big) w(t,x) - f'' \Big( \int_{-\infty}^x v(t,y) \, dy \Big) v(t,x) \big| \, dx
\\
&\le  \Lip f''  \sup_{x\in\RR}\Big\{ \int_{-\infty}^x |w(t,y) - v(t,y)| \,dy \Big\} \|w(t,\cdot) \|_{ L^1(\RR)} \wedge \| v(t,\cdot)\|_{ L^1(\RR)} 
\\
&\qquad+  \| f'' \|_{ L^\infty(\RR)} \int_{\RR} |w(t,x) - v(t,x) | \, dx
\\
&\le  \|w(t) - v(t) \|_{ L^1(\RR)} \big( \Lip f'' \|w_0 \|_{ L^1(\RR)} \wedge \|v_0 \|_{ L^1(\RR)} +  \| f'' \|_{ L^\infty(\RR)} \big).
\endaligned
\]
Also,
\[
\aligned
\Lip k \le \| f'' \|_{ L^\infty(\RR)} \| w \|_{ L^\infty(\RR)} \le  M \| f'' \|_{ L^\infty(\RR)},
\endaligned
\]
where we have used the $L^1$ and $ L^\infty$ bounds \eqref{60},\eqref{40}. This also gives
\[
\aligned
\Theta (\tau, w_0, k) \le \big(\TV(w_0) + \tau \| f''\|_{ L^\infty(\RR)} \| w_0\|_{ L^1(\RR)} \big) e^{M \Lip g \| f''\|_{ L^\infty(\RR)}\tau}
\endaligned
\]
(Recall that $\Theta$ is given by the right-hand side of \eqref{460}). We insert these estimates into \eqref{470} to find
\be
\label{550}
\aligned
\Psi (t) &\le \Psi(0) + a(t) \int_0^t \sup_{[0,\tau)}\Psi(t)\, d\tau,
\endaligned
\ee
for some continuous function $a(t)$ depending on $g,f,w_0$ and $v_0$, which can be explicitly given by the previous estimates.
An application of Gronwall's lemma gives the $L^1$-stability estimate \eqref{50}. This concludes the proof
of Theorem \ref{thm15}. \qed


\section*{Acknowledgements}
We would like to thank Prof.~Jos\'e Francisco Rodrigues for drawing this problem to our attention and for helpful encouragement.
The author was partially supported by the Portuguese Foundation for Science and Technology (FCT)
through the grant PTDC/MAT/110613/2009 - Nonlinear Hyperbolic Systems: Theory and Numerical Approximation, and
by PEst OE/MAT/UI0209/2011.
The author was also supported by the FCT through a 
\emph{Ci\^encia~2008} fellowship. 


\end{document}